\documentclass[12pt,reqno]{article}

\usepackage[usenames]{color}
\usepackage{amssymb}
\usepackage{graphicx}
\usepackage{amscd}

\usepackage[colorlinks=true,
linkcolor=webgreen, filecolor=webbrown,
citecolor=webgreen]{hyperref}

\definecolor{webgreen}{rgb}{0,.5,0}
\definecolor{webbrown}{rgb}{.6,0,0}

\usepackage{color}
\usepackage{float}

\usepackage{graphics,amsmath,amssymb}
\usepackage{amsthm}
\usepackage{amsfonts}
\usepackage{latexsym}
\usepackage{epsf}

\DeclareMathOperator{\ord}{ord}

\begin{document}
\vspace*{5mm}

\begin{center}
\vskip 1cm{\LARGE\bf Additive Combinatorics}\\
\vskip .5cm

{\large\bf with a view towards Computer Science and Cryptography}\\
\vskip .6cm

{\Large\bf An Exposition}\\
\vskip .7cm

\large Khodakhast Bibak\\

\vskip .4cm

Department of Combinatorics and Optimization\\
University of Waterloo\\
Waterloo, Ontario, Canada N2L 3G1\\
\href{mailto: kbibak@uwaterloo.ca}{\tt kbibak@uwaterloo.ca} \vskip
3mm October 25, 2012
\end{center}

\vskip .2 in

\begin{abstract}
Recently, additive combinatorics has blossomed into a vibrant area
in mathematical sciences. But it seems to be a difficult area to
define -- perhaps because of a blend of ideas and techniques from
several seemingly unrelated contexts which are used there. One
might say that additive combinatorics is a branch of mathematics
concerning the study of combinatorial properties of algebraic
objects, for instance, Abelian groups, rings, or fields. This
emerging field has seen tremendous advances over the last few
years, and has recently become a focus of attention among both
mathematicians and computer scientists. This fascinating area has
been enriched by its formidable links to combinatorics, number
theory, harmonic analysis, ergodic theory, and some other
branches; all deeply cross-fertilize each other, holding great
promise for all of them! In this exposition, we attempt to provide
an overview of some breakthroughs in this field, together with a
number of seminal applications to sundry parts of mathematics and
some other disciplines, with emphasis on computer science and
cryptography.
\end{abstract}

\newtheorem{theorem}{Theorem}
\newtheorem{corollary}[theorem]{Corollary}
\newtheorem{lemma}[theorem]{Lemma}
\newtheorem{proposition}[theorem]{Proposition}
\newtheorem{conjecture}[theorem]{Conjecture}
\newtheorem{defin}[theorem]{Definition}
\newenvironment{definition}{\begin{defin}\normalfont\quad}{\end{defin}}
\newtheorem{examp}[theorem]{Example}
\newenvironment{example}{\begin{examp}\normalfont\quad}{\end{examp}}
\newtheorem{rema}[theorem]{Remark}
\newenvironment{remark}{\begin{rema}\normalfont\quad}{\end{rema}}

\newcommand{\len}{\mbox{len}}
\newcommand{\bal}[1]{\begin{align*}#1\end{align*}}
\newcommand{\f}[2]{\displaystyle \frac{#1}{#2}}
\newcommand{\bt}{\begin{thm}}
\newcommand{\et}{\end{thm}}
\newcommand{\bp}{\begin{proof}}
\newcommand{\ep}{\end{proof}}
\newcommand{\bprop}{\begin{prop}}
\newcommand{\eprop}{\end{prop}}
\newcommand{\bl}{\begin{lemma}}
\newcommand{\el}{\end{lemma}}
\newcommand{\bc}{\begin{corollary}}
\newcommand{\ec}{\end{corollary}}
\newcommand{\Z}{\mathbb{Z}}
\newcommand{\be}{\begin{enumerate}}
\newcommand{\ee}{\end{enumerate}}

\section{Introduction}
{\it Additive combinatorics} is a compelling and fast growing area
of research in mathematical sciences, and the goal of this paper
is to survey some of the recent developments and notable
accomplishments of the field, focusing on both pure results and
applications with a view towards computer science and
cryptography. See \cite{TVB} for a book on additive combinatorics,
\cite{NATH1, NATH2} for two books on additive number theory, and
\cite{TRE, Vio1} for two surveys on additive combinatorics. About
additive combinatorics over finite fields and its applications,
the reader is referred to the very recent and excellent survey by
Shparlinski \cite{SH4}.

One might say that additive combinatorics studies combinatorial
properties of algebraic objects, for example, Abelian groups,
rings, or fields, and in fact, focuses on the interplay between
combinatorics, number theory, harmonic analysis, ergodic theory,
and some other branches. Green \cite{GRE} describes additive
combinatorics as the following: ``additive combinatorics is the
study of {\it approximate mathematical structures} such as
approximate groups, rings, fields, polynomials and homomorphisms".
{\it Approximate groups} can be viewed as finite subsets of a
group with the property that they are {\it almost} closed under
multiplication. Approximate groups and their applications (for
example, to expander graphs, group theory, probability, model
theory, and so on) form a very active and promising area of
research in additive combinatorics; the papers \cite{BRGR1, BRGR2,
BGT, BGT2, GRE, HRU, TAO7} contain many developments on this area.
Gowers \cite{GOW4} describes additive combinatorics as the
following: ``additive combinatorics focuses on three classes of
theorems: {\it decomposition theorems}, {\it approximate
structural theorems}, and {\it transference principles}". These
descriptions seem to be mainly inspired by new directions of this
area.

Techniques and approaches applied in additive combinatorics are
often extremely sophisticated, and may have roots in several
unexpected fields of mathematical sciences. For instance,
Hamidoune \cite{HAMI1}, through ideas from connectivity properties
of graphs, established the so-called {\it isoperimetric method},
which is a strong tool in additive combinatorics; see also
\cite{HAMI2, HAMI3, HAMI4} and the nice survey \cite{SERR}. As
another example, Nathanson \cite{NATH3} employed {\it K\"{o}nig's
infinity lemma} on the existence of infinite paths in certain
infinite graphs, and introduced a new class of additive bases, and
also a generalization of the Erd\H{o}s-Tur\'{a}n conjecture about
additive bases of the positive integers. In \cite{PLSC} the
authors employed tools from coding theory to estimating Davenport
constants. Also, in \cite{BABO, MMT1, MMT2, RUZ1, TAO8}
information-theoretic techniques are used to study sumset
inequalities. Very recently, Alon et al. \cite{ABMS1, ABMS2},
using graph-theoretic methods, studied sum-free sets of order $m$
in finite Abelian groups, and also, sum-free subsets of the set
$[1,n]$. Additive combinatorics problems in matrix rings is
another active area of research \cite{BG1, BGS1, BGT, CHAN1,
CHAN5, FHLOSH, GJS, GH1, GH2, HEL1, HEL2, KOW, SV}.

A celebrated result by Szemer\'{e}di, known as Szemer\'{e}di's
theorem (see \cite{AUST1, AUST, BLL, FUR, GOW1, GOW2, GT10, GT1,
NRS1, DHJP, RSC1, RSC2, RSK1, RSK2, SZE2, TAO1, TAO3, TOWS} for
different proofs of this theorem), states that every subset $A$ of
the integers with positive upper density, that is,
$\limsup_{N\rightarrow \infty}|A\cap[1,N]|/N>0$, has arbitrary
long arithmetic progressions. A stunning breakthrough of Green and
Tao \cite{GT2} (that answers a long-standing and folkloric
conjecture by Erd\H{o}s on arithmetic progressions, in a special
case: the primes) says that primes contain arbitrary long
arithmetic progressions. The fusion of methods and ideas from
combinatorics, number theory, harmonic analysis, and ergodic
theory used in its proof is very impressive.

Additive combinatorics has recently found a great deal of
remarkable applications to computer science and cryptography; for
example, to expanders \cite{BKSSW, BRSW, BOU3, BOU12, BG, BG1,
BGS, BGS1, BOVA, BOYE, BGT1, DVSH, DVWI, GUV, KOW, SHEN, VU},
extractors \cite{BIW, BKSSW, BEKO, BEZE, BOU3, BOU7, DGW, DKSM,
DVWI1, GUV, HEHE, LIX, YEH, ZUC1}, pseudorandomness \cite{BV,
LOV1, LMS, LRTV, STV} (also, \cite{TREV0, VAD0} are two surveys
and \cite{VAD1} is a monograph on pseudorandomness), property
testing \cite{BGS0, HASH, HSS, HALO, KALO, KALO1, Sam, SHAP, TUWO}
(see also \cite{GOL}), complexity theory \cite{BLZ, BEZE, BDSS,
BOU5, PLAG1, ZUC1}, hardness amplification \cite{STV, Vio, VW},
probabilistic checkable proofs (PCPs) \cite{ST}, information
theory \cite{BABO, MMT1, MMT2, RUZ1, TAO8}, discrete logarithm
based range protocols \cite{CLS}, non-interactive zero-knowledge
(NIZK) proofs \cite{lip}, compression functions \cite{JOS}, hidden
shifted power problem \cite{BGKS1}, and Diffie-Hellman
distributions \cite{BOU9, BOU1, CFKLLS, FKS}. Additive
combinatorics also has important applications in e-voting
\cite{CLS, lip}. Recently, Bourgain et al. \cite{BDFKK} gave a new
explicit construction of matrices satisfying the {\it Restricted
Isometry Property} (RIP) using ideas from additive combinatorics.
RIP is related to the matrices whose behavior is nearly
orthonormal (at least when acting on sufficiently sparse vectors);
it has several applications, in particular, in compressed sensing
\cite{CAND, CRT, CATA}.

Methods from additive combinatorics provide strong techniques for
studying the so-called {\it threshold phenomena}, which is itself
of significant importance in combinatorics, computer science,
discrete probability, statistical physics, and economics
\cite{ANP, BR, BK, FRI, FRRT, KASA}. There are also very strong
connections between ideas of additive combinatorics and the theory
of {\it random matrices} (see, e.g., \cite{TV1, TV, VU2} and the
references therein); the latter themselves have several
applications in many areas of number theory, combinatorics,
computer science, mathematical and theoretical physics, chemistry,
and so on \cite{ACH, BKSWZ, CT, GMDL, NOVA, RV, TV, TV2}. This
area also has many applications to group theory, analysis,
exponential sums, expanders, complexity theory, discrete geometry,
dynamical systems, and various other scientific disciplines.

Additive combinatorics has seen very fast advancements in the wake
of extremely deep work on Szemer\'{e}di's theorem, the proof of
the existence of long APs in the primes by Green and Tao, and
generalizations and applications of the sum-product problem, and
continues to see significant progress (see \cite{CRLE} for a
collection of open problems in this area). In the next section, we
review Szemer\'{e}di's and Green-Tao theorems (and their
generalizations), two cornerstone breakthroughs in additive
combinatorics. In the third section, we will deal with the
sum-product problem: yet another landmark achievement in additive
combinatorics, and consider its generalizations and applications,
especially to computer science and cryptography.

\section{Szemer\'{e}di's and Green-Tao Theorems, and Their Generalizations}
{\it Ramsey theory} is concerned with the phenomenon that if a
sufficiently large structure (complete graphs, arithmetic
progressions, flat varieties in vector spaces, etc.) is
partitioned arbitrarily into finitely many substructures, then at
least one substructure has necessarily a particular property, and
so total disorder is impossible. In fact, Ramsey theory seeks
general conditions to guarantee the existence of substructures
with regular properties. This theory has many applications, for
example, in number theory, algebra, geometry, topology, functional
analysis (in particular, in Banach space theory), set theory,
logic, ergodic theory, information theory, and theoretical
computer science (see, e.g., \cite{ROS} and the references
therein). {\it Ramsey's theorem} says that in any edge coloring of
a sufficiently large complete graph, one can find monochromatic
complete subgraphs. A nice result of the same spirit is {\it van
der Waerden's theorem}: For a given $k$ and $r$, there exists a
number $N=N(k,r)$ such that if the integers in $[1,N]$ are colored
using $r$ colors, then there is a nontrivial monochromatic
$k$-term arithmetic progression ($k$-AP). Intuitively, this
theorem asserts that in any finite coloring of a structure, one
will find a substructure of the same type at least in one of the
color classes. Note that the finitary and infinitary versions of
the van der Waerden's theorem are equivalent, through a
compactness argument.

One landmark result in Ramsey theory is the Hales-Jewett theorem
\cite{HAJE}, which was initially introduced as a tool for
analyzing certain kinds of games. Before stating this theorem, we
need to define the concept of combinatorial line. A {\it
combinatorial line} is a $k$-subset in the $n$-dimensional grid
$[1,k]^{n}$ yielded from some template in $([1,k]\cup \{*\})^{n}$
by replacing the symbol $*$ with $1, \ldots, k$ in turn. The
Hales-Jewett theorem states that for every $r$ and $k$ there
exists $n$ such that every $r$-coloring of the $n$-dimensional
grid $[1,k]^{n}$ contains a combinatorial line. Roughly speaking,
it says that for every multidimensional grid whose faces are
colored with a number of colors, there must necessarily be a line
of faces of all the same color, if the dimension is sufficiently
large (depending on the number of sides of the grid and the number
of colors). Note that instead of seeking arithmetic progressions,
the Hales-Jewett theorem seeks combinatorial lines. This theorem
has many interesting consequences in Ramsey theory, two of which,
are van der Waerden's theorem and its multidimensional version,
i.e., the Gallai-Witt theorem (see, e.g., \cite{GKP, GRS, JUK} for
further information).

Erd\H{o}s and Tur\'{a}n \cite{ET} proposed a very strong form of
van der Waerden's theorem -- the density version of van der
Waerden's theorem. They conjectured that arbitrarily long APs
appear not only in finite partitions but also in every
sufficiently dense subset of positive integers. More precisely,
the Erd\H{o}s-Tur\'{a}n conjecture states that if $\delta$ and $k$
are given, then there is a number $N=N(k,\delta)$ such that any
set $A \subseteq [1,N]$ with $|A|\geq \delta N$ contains a
non-trivial $k$-AP. Roth \cite{ROT} employed methods from Fourier
analysis (or more specifically, the {\it Hardy-Littlewood circle
method}) to prove the $k=3$ case of the Erd\H{o}s-Tur\'{a}n
conjecture (see also \cite{BLOO, BOU2, BOU10, CRSI, ELK, GRWO,
LEV, MESH, OBR, SAN1, SAN2, SHKR1, SOL6}). Szemer\'{e}di
\cite{SZE1} verified the Erd\H{o}s-Tur\'{a}n conjecture for
arithmetic progressions of length four. Finally, Szemer\'{e}di
\cite{SZE2} by a tour de force of sophisticated combinatorial
arguments proved the conjecture, now known as {\it Szemer\'{e}di's
theorem} -- one of the milestones of combinatorics. Roughly
speaking, this theorem states that long arithmetic progressions
are very widespread and in fact it is not possible to completely
get rid of them from a set of positive integers unless we can
contract the set (sufficiently) to make it of density zero.
{\L}aba and Pramanik \cite{LAPR} (also see \cite{POT}) proved that
every compact set of reals with Lebesgue measure zero supporting a
probabilistic measure satisfying appropriate dimensionality and
Fourier decay conditions must contain non-trivial 3APs.

Conlon and Gowers \cite{COGO} considered Szemer\'{e}di's theorem,
and also several other combinatorial theorems such as Tur\'{a}n's
theorem and Ramsey's theorem in sparse random sets. Also,
Szemer\'{e}di-type problems in various structures other than
integers have been a focus of significant amount of work. For
instance, \cite{GRE1, LIWO} consider these kinds of problems in
the finite field setting. Very recently, Bateman and Katz
\cite{BAKA1} (also see \cite{BAKA2}) achieved new bounds for the
{\it cap set problem}, which is basically Roth's problem, but in a
vector space over finite fields (a set $A\subset \mathbb{F}_3^N$
is called a {\it cap set} if it contains no lines).

A salient ingredient in Szemer\'{e}di's proof (in addition to van
der Waerden's theorem) is the {\it Szemer\'{e}di regularity
lemma}. This lemma was conceived specifically for the purpose of
this proof, but is now, by itself, one of the most powerful tools
in extremal graph theory (see, e.g., \cite{KS, VRMS2}, which are
two surveys on this lemma and its applications). Roughly speaking,
it asserts that the vertex set of every (large) graph can be
partitioned into relatively few parts such that the subgraphs
between the parts are random-like. Indeed, this result states that
each large dense graph may be decomposed into a low-complexity
part and a pseudorandom part (note that Szemer\'{e}di's regularity
lemma is the archetypal example of the {\it dichotomy between
structure and randomness} \cite{TAO6}). The lemma has found
numerous applications not only in graph theory, but also in
discrete geometry, additive combinatorics, and computer science.
For example, as Trevisan \cite{TRE} mentions, to solve a
computational problem on a given graph, it might be easier to
first construct a Szemer\'{e}di approximation -- this resulted
approximating graph has a simpler configuration and would be
easier to treat. Note that the significance of Szemer\'{e}di's
regularity lemma goes beyond graph theory: it can be reformulated
as a result in information theory, approximation theory, as a
compactness result on the completion of the space of finite
graphs, etc. (see \cite{LLS} and the references therein). Very
recently, Tao and Green \cite{GT10} established an {\it arithmetic
regularity lemma} and a complementary {\it arithmetic counting
lemma} that have several applications, in particular, an
astonishing proof of Szemer\'{e}di's theorem.

The {\it triangle removal lemma} established by Ruzsa and
Szemer\'{e}di \cite{IZRES} is one of the most notable applications
of Szemer\'{e}di's regularity lemma. It asserts that each graph of
order $n$ with $o(n^{3})$ triangles can be made triangle-free by
removing $o(n^{2})$ edges. In other words, if a graph has
asymptotically few triangles then it is asymptotically close to
being triangle-free. As a clever application of this lemma, Ruzsa
and Szemer\'{e}di \cite{IZRES} obtained a new proof of Roth's
theorem (see also \cite{SOL6}, in which the author using the
triangle removal lemma proves Roth type theorems in finite
groups). Note that a generalization of the triangle removal lemma,
known as {\it simplex removal lemma}, can be used to deduce
Szemer\'{e}di's theorem (see \cite{GOW2, RSK1, RSK2, TAO3}). The
triangle removal lemma was extended by Erd\H{o}s, Frankl, and
R\"{o}dl \cite{EFR} to the {\it graph removal lemma}, which
roughly speaking, asserts that if a given graph does not contain
too many subgraphs of a given type, then all the subgraphs of this
type can be removed by deleting a few edges. More precisely, given
a fixed graph $H$ of order $k$, any graph of order $n$ with
$o(n^{k})$ copies of $H$ can be made $H$-free by removing
$o(n^{2})$ edges. Fox \cite{FOX} gave a proof of the graph removal
lemma which avoids applying Szemer\'{e}di's regularity lemma and
gives a better bound (also see \cite{COFO}). The graph removal
lemma has many applications in graph theory, additive
combinatorics, discrete geometry, and theoretical computer
science. One surprising application of this lemma is to the area
of {\it property testing}, which is now a very dynamic area in
computer science \cite{AFNS, AS0, AS1, AS2, AUTA, GOL, GGR, HKMS,
VRMS, VRMS1, RON, RON0, Sam}. Property testing typically refers to
the existence of sub-linear time probabilistic algorithms (called
testers), which distinguish between objects $G$ (e.g., a graph)
having a given property $P$ (e.g., bipartiteness) and those being
far away (in an appropriate metric) from $P$. Property testing
algorithms have been recently designed and utilized for many kinds
of objects and properties, in particular, discrete properties
(e.g., graph properties, discrete functions, and sets of
integers), geometric properties, algebraic properties, etc.

There is a (growing) number of proofs of Szemer\'{e}di's theorem,
arguably seventeen proofs to this date. One such elegant proof
that uses ideas from model theory was given by Towsner
\cite{TOWS}. For another model theory based perspective, see
\cite{MASH}, in which the authors give stronger regularity lemmas
for some classes of graphs.

In fact, one might claim that many of these proofs have themselves
opened up a new field of research. Furstenberg \cite{FUR} by
rephrasing it as a problem in {\it dynamical systems}, and then
applying several powerful techniques from ergodic theory achieved
a nice proof of the Szemer\'{e}di's theorem. In fact, Furstenberg
presented a correspondence between problems in the subsets of
positive density in the integers and recurrence problems for sets
of positive measure in a {\it probability measure preserving
system}. This observation is now known as the {\it Furstenberg
correspondence principle}. {\it Ergodic theory} is concerned with
the long-term behavior in dynamical systems from a statistical
point of view (see, e.g., \cite{EIWA}). This area and its
formidable way of thinking have made many strong connections with
several branches of mathematics, including combinatorics, number
theory, coding theory, group theory, and harmonic analysis; see,
for example, \cite{KRA, KRA1, KRA3, KRA2} and the references
therein for some connections between ergodic theory and additive
combinatorics.

This ergodic-theoretic method is one of the most flexible known
proofs, and has been very successful at reaching considerable
generalizations of Szemer\'{e}di's theorem. Furstenberg and
Katznelson \cite{FK} obtained the {\it multidimensional
Szemer\'{e}di theorem}. Their proof relies on the concept of {\it
multiple recurrence}, a powerful tool in the interaction between
ergodic theory and additive combinatorics. A purely combinatorial
proof of this theorem was obtained roughly in parallel by Gowers
\cite{GOW2}, and Nagle et al. \cite{NRS1, RSC1, RSC2, RSK1, RSK2},
and subsequently by Tao \cite{TAO3}, via establishing a {\it
hypergraph removal lemma} (see also \cite{VRMS1, TAO4}). Also,
Austin \cite{AUST1} proved the theorem via both ergodic-theoretic
and combinatorial approaches. The multidimensional Szemer\'{e}di
theorem was significantly generalized by Furstenberg and
Katznelson \cite{FK2} (via ergodic-theoretic approaches), and
Austin \cite{AUST} (via both ergodic-theoretic and combinatorial
approaches), to the {\it density Hales-Jewett theorem}. The
density Hales-Jewett theorem states that for every $\delta > 0$
there is some $N_{0} \geq 1$ such that whenever $A \subseteq
[1,k]^N$ with $N \geq N_{0}$ and $|A|\geq \delta k^N$, $A$
contains a combinatorial line. Recently, in a massively
collaborative online project, namely {\it Polymath 1} (a project
that originated in Gowers' blog), the Polymath team found a purely
combinatorial proof of the density Hales-Jewett theorem, which is
also the first one providing explicit bounds for how large $n$
needs to be \cite{DHJP} (also see \cite{NIE}). Such bounds could
not be obtained through the ergodic-theoretic methods, since these
proofs rely on the Axiom of Choice. It is worth mentioning that
this project was selected as one of the TIME Magazine's Best Ideas
of 2009.

Furstenberg's proof gave rise to the field of {\it ergodic Ramsey
theory}, in which arithmetical, combinatorial, and geometrical
configurations preserved in (sufficiently large) substructures of
a structure, are treated via ideas and techniques from ergodic
theory (or more specifically, multiple recurrence). Ergodic Ramsey
theory has since produced a high number of combinatorial results,
some of which have yet to be obtained by other means, and has also
given a deeper understanding of the structure of measure
preserving systems. In fact, ergodic theory has been used to solve
problems in Ramsey theory, and reciprocally, Ramsey theory has led
to the discovery of new phenomena in ergodic theory. However, the
ergodic-theoretic methods and the infinitary nature of their
techniques have some limitations. For example, these methods do
not provide any effective bound, since, as we already mentioned,
they rely on the Axiom of Choice. Also, despite van der Waerden's
theorem is not directly used in Furstenberg's proof, probably any
effort to make the proof quantitative would result in rapidly
growing functions. Furthermore, the ergodic-theoretic methods, to
this day, have the limitation of only being able to deal with sets
of positive density in the integers, although this density is
allowed to be arbitrarily small. However, Green and Tao \cite{GT2}
discovered a {\it transference principle} which allowed one to
reduce problems on structures in special sets of zero density
(such as the primes) to problems on sets of positive density in
the integers. It is worth mentioning that Conlon, Fox, and Zhao
\cite{CFZ} established a {\it transference principle} extending
several classical extremal graph theoretic results, including the
removal lemmas for graphs and groups (the latter leads to an
extension of Roth's theorem), the Erd\H{o}s-Stone-Simonovits
theorem and Ramsey's theorem, to sparse pseudorandom graphs.

Gowers \cite{GOW1} generalized the arguments previously studied in
\cite{GOW3, ROT}, in a substantial way. In fact, he employed
combinatorics, generalized Fourier analysis, and inverse
arithmetic combinatorics (including multilinear versions of {\it
Freiman's theorem} on sumsets, and the {\it Balog-Szemer\'{e}di
theorem}) to reprove Szemer\'{e}di's theorem with explicit bounds.
Note that Fourier analysis has a wide range of applications, in
particular, to cryptography, hardness of approximation, signal
processing, threshold phenomena for probabilistic models such as
random graphs and percolations, and many other disciplines.
Gowers' article introduced a kind of {\it higher degree Fourier
analysis}, which has been further developed by Green and Tao.
Indeed, Gowers initiated the study of a new measure of functions,
now referred to as {\it Gowers (uniformity) norms}, that resulted
in a better understanding of the notion of {\it pseudorandomness}.

The Gowers norm, which is an important special case of noise
correlation (intuitively, the {\it noise correlation} between two
functions $f$ and $g$ measures how much $f(x)$ and $g(y)$
correlate on random inputs $x$ and $y$ which are correlated),
enjoys many properties and applications, and is now a very dynamic
area of research in mathematical sciences; see \cite{AUM, GT5,
GT4, GT8, GTZ1, GTZ2, HASH, HALO, KALO, LIU, LOV, LMS0, MATT, Sam,
ST} for more properties and applications of the Gowers norm. Also,
the best known bounds for Szemer\'{e}di's theorem are obtained
through the so-called {\it inverse theorems} for Gowers norms.
Recently, Green and Tao \cite{GT1} (see also \cite{TAO10}), using
the density-increment strategy of Roth \cite{ROT} and Gowers
\cite{GOW3, GOW1}, derived Szemer\'{e}di's theorem from the {\it
inverse conjectures GI($s$)} for the Gowers norms, which were
recently established in \cite{GTZ2}.

To the best of my knowledge, there are two types of inverse
theorems in additive combinatorics, namely the {\it inverse sumset
theorems of Freiman type} (see, e.g., \cite{CHAN4, FGSS1, FGSS2,
GRRU, GT12, SCHO, SSV, SZVU1, SZVU2, SZVU3, TAO9} and \cite{FREI,
NATH2}), and {\it inverse theorems for the Gowers norms} (see,
e.g., \cite{GOW1, GOWO1, GOWO2, GOWO3, GT5, GT7, GT6, GT8, GTZ1,
GTZ2, HASH, HALO, HOKR1, KALO, LOV, LMS0, MATT, Sam, TZ2, TUWO}).
It is interesting that the inverse conjecture leads to a finite
field version of Szemer\'{e}di's theorem \cite{TAO10}: Let
$\mathbb{F}_p$ be a finite field. Suppose that $\delta > 0$, and
$A \subset \mathbb{F}_p^n$ with $|A|\geq \delta |\mathbb{F}_p^n|$.
If $n$ is sufficiently large depending on $p$ and $\delta$, then
$A$ contains an (affine) line $\{x, x + r, \ldots, x + (p - 1)r\}$
for some $x, r \in \mathbb{F}_p^n$ with $r\not=0$ (actually, $A$
contains an affine $k$-dimensional subspace, $k\geq 1$).

Suppose $r_k(N)$ is the cardinality of the largest subset of
$[1,N]$ containing no nontrivial $k$-APs. Giving asymptotic
estimates on $r_k(N)$ is an important inverse problem in additive
combinatorics. Behrend \cite{BEH} proved that
$$r_3(N)=\Omega \Bigg(\frac{N}{2^{2\sqrt{2}\sqrt{\log_2 N}}.\log^{1/4}
N}\Bigg).$$ Rankin \cite{RAN} generalized Behrend's construction
to longer APs. Roth proved that $r_3(N)=o(N)$. In fact, he proved
the first nontrivial upper bound
$$r_3(N)=O\Bigg(\frac{N}{\log\log N}\Bigg).$$ Bourgain \cite{BOU2,
BOU10} improved Roth's bound. In fact, Bourgain \cite{BOU10} gave
the upper bound
$$r_3(N)=O\Bigg(\frac{N(\log\log N)^2}{\log^{2/3}N}\Bigg).$$
Sanders \cite{SAN2} proved the following upper bound which is the
state-of-the-art:
$$r_3(N)=O\Bigg(\frac{N(\log\log N)^5}{\log N}\Bigg).$$
Bloom \cite{BLOO} through the nice technique ``translation of a
proof in $\mathbb{F}_q[t]$ to one in $\mathbb{Z}/N\mathbb{Z}$",
extends Sanders' proof to 4 and 5 variables. As Bloom mentions in
his paper, many problems of additive combinatorics might be easier
to attack via the approach ``translating from $\mathbb{F}_p^N$ to
$\mathbb{F}_q[t]$ and hence to $\mathbb{Z}/N\mathbb{Z}$".

Elkin \cite{ELK} managed to improve Behrend's 62-year old lower
bound by a factor of $\Theta(\log^{1/2}N)$. Actually, Elkin showed
that
$$r_3(N)=\Omega \Bigg(\frac{N}{2^{2\sqrt{2}\sqrt{\log_2 N}}}.\log^{1/4} N \Bigg).$$
See also \cite{GRWO} for a short proof of Elkin's result, and
\cite{OBR} for constructive lower bounds for $r_k(N)$. Schoen and
Shkredov \cite{SCSH} using ideas from the paper of Sanders
\cite{SAN3} and also the new probabilistic technique established
by Croot and Sisask \cite{CRSI}, obtained Behrend-type bounds for
linear equations involving 6 or more variables. Thanks to this
result, one may see that perhaps the Behrend-type constructions
are not too far from being best-possible.

Almost all the known proofs of Szemer\'{e}di's theorem are based
on a {\it dichotomy between structure and randomness} \cite{TAO6,
TAO5}, which allows many mathematical objects to be split into a
`structured part' (or `low-complexity part') and a `random part'
(or `discorrelated part'). Tao \cite{TAO1} best describes almost
all known proofs of Szemer\'{e}di's theorem collectively as the
following: ``Start with the set $A$ (or some other object which is
a proxy for $A$, e.g., a graph, a hypergraph, or a
measure-preserving system). For the object under consideration,
define some concept of randomness (e.g., $\varepsilon$-regularity,
uniformity, small Fourier coefficients, or weak mixing), and some
concept of structure (e.g., a nested sequence of arithmetically
structured sets such as progressions or Bohr sets, or a partition
of a vertex set into a controlled number of pieces, a collection
of large Fourier coefficients, a sequence of almost periodic
functions, a tower of compact extensions of the trivial factors).
Obtain some sort of structure theorem that splits the object into
a structured component, plus an error which is random relative to
that structured component. To prove Szemer\'{e}di's theorem (or a
variant thereof), one then needs to obtain some sort of {\it
generalized von Neumann theorem} \cite{GT2} to eliminate the
random error, and then some sort of {\it structured recurrence
theorem} for the structured component".

Erd\H{o}s's famous conjecture on APs states that a set $A = \{a_1,
a_2, \ldots, a_n, \ldots\}$ of positive integers, where $a_i <
a_{i+1}$ for all $i$, with the divergent sum $\sum_{n\in
\mathbb{Z}^{+}} \frac{1}{a_n}$, contains arbitrarily long APs. If
true, the theorem includes both Szemer\'{e}di's and Green-Tao
theorems as special cases. This conjecture seems to be too strong
to hold, and in fact, might be very difficult to attack -- it is
not even known whether such a set must contain a 3-AP! So, let us
mention an equivalent statement for Erd\H{o}s's conjecture that
may be helpful. Let $N$ be a positive integer. For a positive
integer $k$, define $a_{k}(N):=r_{k}(N)/ N$ (note that
Szemer\'{e}di's theorem asserts that $\lim_{N\rightarrow \infty}
a_{k}(N)=0$, for all $k$). It can be proved (see \cite{SHKR}) that
Erd\H{o}s's conjecture is true if and only if the series
$\sum_{i=1}^{\infty} a_{k}4^{i}$ converges for any integer $k\geq
3$. So, to prove Erd\H{o}s's conjecture, it suffices to obtain the
estimate $a_{k}(N)\ll 1/(\log N)^{1+\varepsilon}$, for any $k\geq
3$ and for some $\varepsilon >0$.

Szemer\'{e}di's theorem plays an important role in the proof of
the Green-Tao theorem \cite{GT2}: The primes contain arithmetic
progressions of arbitrarily large length (note that the same
result is valid for every subset of the primes with positive
relative upper density). Green and Tao \cite{GT3} also proved that
there is a $k$-AP of primes all of whose terms are bounded by
$$2^{2^{2^{2^{2^{2^{2^{2^{(100k)}}}}}}}},$$ which shows that how
far out in the primes one must go to warrant finding a $k$-AP. A
conjecture (see \cite{KRA}) asserts that there is a $k$-AP in the
primes all of whose terms are bounded by $k!+1$.

There are three fundamental ingredients in the proof of the
Green-Tao theorem (in fact, there are many similarities between
Green and Tao's approach and the ergodic-theoretic method, see
\cite{HOKR}). The first is Szemer\'{e}di's theorem itself. Since
the primes do not have positive upper density, Szemer\'{e}di's
theorem cannot be directly applied. The second major ingredient in
the proof is a certain {\it transference principle} that allows
one to use Szemer\'{e}di's theorem in a more general setting (a
generalization of Szemer\'{e}di's theorem to the {\it pseudorandom
sets}, which can have zero density). The last major ingredient is
applying some notable features of the primes and their
distribution through results of Goldston and Yildirim \cite{GY1,
GY2}, and proving the fact that this generalized Szemer\'{e}di
theorem can be efficiently applied to the primes, and indeed, the
set of primes will have the desired pseudorandom properties.

In fact, Green and Tao's proof employs the techniques applied in
several known proofs of Szemer\'{e}di's theorem and exploits a
dichotomy between structure and randomness. This proof is based on
ideas and results from several branches of mathematics, for
example, combinatorics, analytical number theory,
pseudorandomness, harmonic analysis, and ergodic theory. Reingold
et al. \cite{RTTV}, and Gowers \cite{GOW4}, independently obtained
a short proof for a fundamental ingredient of this proof.

Tao and Ziegler \cite{TZ} (see also \cite{RTTV}), via a
transference principle for polynomial configurations, extended the
Green-Tao theorem to cover polynomial progressions: Let $A\subset
{\cal P}$ be a set of primes of positive relative upper density in
the primes, i.e., $\limsup_{N\rightarrow \infty} |A\cap
[1,N]|/|{\cal P}\cap [1,N]|>0$. Then, given any integer-valued
polynomials $P_1,\ldots,P_k$ in one unknown $m$ with vanishing
constant terms, the set $A$ contains infinitely many progressions
of the form $x+P_1(m),\ldots, x+P_k(m)$ with $m>0$ (note that the
special case when the polynomials are $m, 2m,\ldots, km$ implies
the previous result that there are $k$-APs of primes). Tao
\cite{TAO2} proved the analogue in the Gaussian integers. Green
and Tao (in view of the parallelism between the integers and the
polynomials over a finite field) thought that the analogue of
their theorem should be held in the setting of function fields; a
result that was proved by L\^{e} \cite{LE}: Let $\mathbb{F}_{q}$
be a finite field over $q$ elements. Then for any $k>0$, one can
find polynomials $f,g \in \mathbb{F}_{q}[t], \;g \not \equiv 0$
such that the polynomials $f +Pg$ are all irreducible, where $P$
runs over all polynomials $P \in \mathbb{F}_{q}[t]$ of degree less
than $k$. Moreover, such structures can be found in every set of
positive relative upper density among the irreducible polynomials.
The proof of this interesting theorem follows the ideas of the
proof of the Green-Tao theorem very closely.

\section{Sum-Product Problem: Its Generalizations and Applications}
The {\it sum-product problem} and its generalizations constitute
another vibrant area in additive combinatorics, and have led to
many seminal applications to number theory, Ramsey theory,
computer science, and cryptography.

Let's start with the definition of sumset, product set, and some
preliminaries. We will follow closely the presentation of Tao
\cite{TAO5}. Let $A$ be a finite nonempty set of elements of a
ring $R$. We define the {\it sumset} $A+A=\{a+b\; :\; a,b\in A
\},$ and the {\it product set} $A\cdot A=\{a\cdot b \;:\; a,b\in A
\}$. Suppose that no $a\in A$ is a zero divisor (otherwise,
$A\cdot A$ may become very small, which lead to degenerate cases).
Then one can easily show that $A+A$ and $A\cdot A$ will be at
least as large as $A$. The set $A$ may be almost closed under
addition, which, for example, occurs when $A$ is an arithmetic
progression or an additive subgroup in the ring $R$ (e.g., if
$A\subset \mathbb{R}$ is an AP, then $|A+A|=2|A|-1$, and $|A\cdot
A|\geq c|A|^{2-\varepsilon}$), or it may be almost closed under
multiplication, which, for example, occurs when $A$ is a geometric
progression or a multiplicative subgroup in the ring $R$ (e.g., if
$N\subset \mathbb{R}$ is an AP and $A = \{2^{n} \; :\; n \in N\}$,
then $|A\cdot A|=2|A|-1$, and $|A+A|\approx |A|^{2}$). Note that
even if $A$ is a dense subset of an arithmetic progression or
additive subgroup (or a dense subset of an geometric progression
or multiplicative subgroup), then $A+A$ (or $A\cdot A$,
respectively) is still comparable in size to $A$. But it is
difficult for $A$ to be almost closed under addition and
multiplication simultaneously, unless it is very close to a
subring. The {\it sum-product phenomenon} says that if a finite
set $A$ is not close to a subring, then either the sumset $A+A$ or
the product set $A\cdot A$ must be considerably larger than $A$.
The reader can refer to \cite{SCSH1} and the references therein to
see some lower bounds on $|C-C|$ and $|C+C|$, where $C$ is a
convex set (a set of reals $C=\{c_1, \ldots , c_n\}$ is called
{\it convex} if $c_{i+1}-c_i>c_i-c_{i-1}$, for all $i$).

In the reals setting, does there exist an $A\subset \mathbb{R}$
for which $\max \{|A+A|,|A\cdot A|\}$ is `small'? Erd\H{o}s and
Szemer\'{e}di \cite{ES} gave a negative answer to this question.
Actually, they proved the inequality $\max \{|A+A|,|A\cdot
A|\}\geq c|A|^{1+\varepsilon}$ for a small but positive
$\varepsilon$, where $A$ is a subset of the reals. They also
conjectured that $\max \{|A+A|,|A\cdot A|\}\geq c|A|^{2-\delta}$,
for any positive $\delta$. Much efforts have been made towards the
value of $\varepsilon$. Elekes \cite{ELE} observed that the
sum-product problem has interesting connections to problems in
incidence geometry. In particular, he applied the so-called {\it
Szemer\'{e}di-Trotter theorem} and showed that $\varepsilon \geq
1/4$, if $A$ is a finite set of real numbers. Elekes's result was
extended to complex numbers in \cite{SOL3}. In the case of reals,
the state-of-the-art is due to Solymosi \cite{SOL4}: one can take
$\varepsilon$ arbitrarily close to $1/3$. For complex numbers,
Solymosi \cite{SOL1}, using the Szemer\'{e}di-Trotter theorem,
proved that one can take $\varepsilon$ arbitrarily close to
$3/11$. Very recently, Rudnev \cite{RUDN1}, again using the
Szemer\'{e}di-Trotter theorem, obtained a bound which is the
state-of-the-art in the case of complex numbers: one can take
$\varepsilon$ arbitrarily close to $19/69$.

Solymosi and Vu \cite{SV} proved a sum-product estimate for a
special finite set of square matrices with complex entries, where
that set is {\it well-conditioned} (that is, its matrices are far
from being singular). Note that If we remove the latter condition
(i.e., well-conditioned!) then the theorem will not be true; see
\cite[Example 1.1]{SV}.

Wolff \cite{WOL} motivated by the `finite field Kakeya
conjecture', formulated the finite field version of sum-product
problem. The {\it Kakeya conjecture} says that the Hausdorff
dimension of any subset of $\mathbb{R}^n$ that contains a unit
line segment in every direction is equal to $n$; it is open in
dimensions at least three. The {\it finite field Kakeya
conjecture} asks for the smallest subset of $\mathbb{F}_q^n$ that
contains a line in each direction. This conjecture was proved by
Dvir \cite{DVIR0} using a clever application of the so-called {\it
polynomial method}; see also \cite{DVIR} for a nice survey on this
problem and its applications especially in the area of randomness
extractors. The polynomial method, which has proved to be very
useful in additive combinatorics, is roughly described as the
following: Given a field $\mathbb{F}$ and a finite subset
$S\subset \mathbb{F}^n$. Multivariate polynomials over
$\mathbb{F}$ which vanish on all points of $S$, usually get some
combinatorial properties about $S$. (This has some similarities
with what we usually do in algebraic geometry!) See, e.g.,
\cite[Chapter 16]{JUK} for some basic facts about the polynomial
method, \cite{DVIR0, GUT, GUKA1, GUKA2} for applications in
additive combinatorics, and \cite{DVWI1, GUV, SUDA} for
applications in computer science.

Actually, the finite field version (of sum-product problem)
becomes more difficult, because we will encounter with some
difficulties in applying the Szemer\'{e}di-Trotter incidence
theorem in this setting. In fact, the {\it crossing lemma}, which
is an important ingredient in the proof of Szemer\'{e}di-Trotter
theorem \cite{LSY0}, relies on {\it Euler's formula} (and so on
the topology of the plane), and consequently does not work in
finite fields. Note that the proof that Szemer\'{e}di and Trotter
presented for their theorem was somewhat complicated, using a
combinatorial technique known as {\it cell decomposition}
\cite{STJR}.

When working with finite fields it is important to consider fields
whose order is prime and not the power of a prime; because in the
latter case we can take $A$ to be a subring which leads to the
degenerate case $|A|=|A+A|=|A\cdot A|$. A stunning result in the
case of finite field $\mathbb{F}_{p}$, with $p$ prime, was proved
by Bourgain, Katz and Tao \cite{BKT}. They proved the following:

if $A\subset \mathbb{F}_{p}$, and $p^{\delta}\leq |A| \leq
p^{1-\delta}$ for some $\delta > 0$, then there exists
$\varepsilon=\varepsilon(\delta)> 0$ such that $\max
\{|A+A|,|A\cdot A|\}\geq c|A|^{1+\varepsilon}$.

This result is now known as the {\it sum-product theorem} for
$\mathbb{F}_{p}$. In fact, this theorem holds if $A$ is not too
close to be the whole field. The condition $|A| \geq p^{\delta}$
in this theorem was removed by Bourgain, Glibichuk, and Konyagin
in \cite{BGK}. Also, note that the condition $|A| \leq
p^{1-\delta}$ is necessary (e.g., if we consider a set $A$
consisting of all elements of the field except one, then $\max
\{|A+A|,|A\cdot A|\}=|A|+1$). The idea for the proof of this
theorem is by contradiction; assume that $|A+A|$ and $|A\cdot A|$
are close to $|A|$ and conclude that $A$ is behaving very much
like a subfield of $\mathbb{F}_{p}$. Sum-product estimates for
rational functions (i.e., the results that one of $A+A$ or $f(A)$
is substantially larger than $A$, where $f$ is a rational
function) have also been treated (see, e.g., \cite{BIW, BT}).
Also, note that problems of the kind `interaction of summation and
addition' are very important in various contexts of additive
combinatorics and have very interesting applications (see
\cite{BAL, GASH, GLI, GLRU, GUKA2, IRR, OSSH, RORU}).

Garaev \cite{GAR} proved the first quantitative sum-product
estimate for fields of prime order: Let $A\subset \mathbb{F}_p$
such that $1<|A|<p^{7/13}\log^{-4/13}p$. Then $$\max
\{|A+A|,|A\cdot A|\}\gg \frac{|A|^{15/14}}{\log^{2/7}|A|}.$$
Garaev's result was extended and improved by several authors.
Rudnev \cite{RUDN} proved the following: Let $A\subset
\mathbb{F}_p^*$ with $|A|<\sqrt{p}$ and $p$ large. Then $$\max
\{|A+A|,|A\cdot A|\}\gg \frac{|A|^{12/11}}{\log ^{4/11}|A|}.$$ Li
and Roche-Newton \cite{LIRO} proved a sum-product estimate for
subsets of a finite field whose order is not prime: Let $A\subset
\mathbb{F}_{p^{n}}$ with $|A \cap cG|\leq |G|^{1/2}$ for any
subfield $G$ of $\mathbb{F}_{p^{n}}$ and any element $c \in
\mathbb{F}_{p^{n}}$. Then $$\max \{|A+A|,|A\cdot A|\}\gg
\frac{|A|^{12/11}}{\log_2 ^{5/11}|A|}.$$ See also \cite{BOGAR,
GAR1, GAR3, KASH1, KASH2, LI, SHEN1, SHEN2} for other
generalizations and improvements of Garaev's result. As an
application, Shparlinski \cite{SH3} using Rudnev's result
\cite{RUDN}, estimates the cardinality, $\#\Gamma_p(T)$, of the
set
$$\Gamma_p(T)=\{\gamma \in \mathbb{F}_p\; :
\;  \ord\gamma \leq T \; \text{and} \;
\ord(\gamma+\gamma^{-1})\leq T\},$$ where $\ord\gamma$
(multiplicative order of $\gamma$) is the smallest positive
integer $t$ with $\gamma^t=1$.

T\'{o}th \cite{TOT} generalized the Szemer\'{e}di-Trotter theorem
to complex points and lines in $\mathbb{C}^2$ (also see \cite{ZAH}
for a different proof of this result, and a sharp result in the
case of $\mathbb{R}^4$). As another application of the sum-product
theorem, Bourgain, Katz and Tao \cite{BKT} (also see \cite{BOU3})
derived an important Szemer\'{e}di-Trotter type theorem in prime
finite fields (but did not quantify it):

If $\mathbb{F}_p$ is a prime field, and ${\cal P}$ and ${\cal L}$
are points and lines in the projective plane over $\mathbb{F}_p$
with cardinality $|{\cal P}|, |{\cal L}| \leq N < p^{\alpha}$ for
some $0 < \alpha < 2$, then $\bigl|\{(p,l)\in {\cal P} \times
{\cal L} \;:\; p \in l\}\bigr|\leq CN^{3/2-\varepsilon}$, for some
$\varepsilon=\varepsilon(\alpha)>0$.

Note that it is not difficult to generalize this theorem from
prime finite fields to every finite field that does not contain a
large subfield. The first quantitative Szemer\'{e}di-Trotter type
theorem in prime finite fields was obtained by Helfgott and Rudnev
\cite{HERU}. They showed that $\varepsilon \geq 1/10678$, when
$|{\cal P}|=|{\cal L}|<p$. (Note that this condition prevents $P$
from being the entire plane $\mathbb{F}_p^2$.) This result was
extended to general finite fields (with a slightly weaker
exponent) by Jones \cite{JONE}. Also, Jones \cite{JONE1, JONE3}
improved the result of Helfgott and Rudnev \cite{HERU} by
replacing $1/10678$ with $1/806-o(1)$, and $1/662-o(1)$,
respectively. A near-sharp generalization of the
Szemer\'{e}di-Trotter theorem to higher dimensional points and
varieties was obtained in \cite{SOTA}.

Using ideas from additive combinatorics (in fact, combining the
techniques of cell decomposition and polynomial method in a novel
way), Guth and Katz \cite{GUKA2}, achieved a near-optimal bound
for the {\it Erd\H{o}s distinct distance problem} in the plane.
They proved that a set of $N$ points in the plane has at least
$c{N \over \log N}$ distinct distances (see also \cite{GIS} for
some techniques and ideas related to this problem).

Hart, Iosevich, and Solymosi \cite{HIS} obtained a new proof of
the sum-product theorem based on incidence theorems for hyperbolas
in finite fields which is achieved through some estimates on
Kloosterman exponential sums. Some other results related to the
incidence theorems can be found in \cite{CHIKR, HIKR}.

The sum-product theorem has a plethora of deep applications to
various areas such as incidence geometry \cite{BDSS, BOU3, BKT,
CHSO, CEHIK, IRR, JONE2, RORU, TAO0, VU, VWW}, analysis
\cite{CHAN, GJS, HLS, TAO0}, PDE \cite{TAO0}, group theory
\cite{BOVA, BGT, GH1, GH2, HEL1, HEL2, HRU, KOW}, exponential and
character sums \cite{BBS1, BOU9, BOU1, BOU3, BOU5, BOU6, BOU8,
BOU14, BOU4, BOU11, BOCH1, BGKS1, BGKS2, BGK, CHAN3, GAR2}, number
theory \cite{BBS, BBS1, BOU8, BOU4, BGS, BGS1, CHAN3, CCGHSZ,
CIGA, CGOS, CSZ, OSSH}, combinatorics \cite{BKSSW, BOVA, BG1,
TAO0}, expanders \cite{BOU3, BG, BG1, BGS, BGS1, BOVA, BGT1, KOW,
SHEN, VU}, extractors \cite{BIW, BKSSW, BOU3, BOU7, CG, DGW, HEHE,
YEH}, dispersers \cite{BKSSW}, complexity theory \cite{BOU5},
pseudorandomness \cite{BV, BOU3, LRTV}, property testing
\cite{GOL, Sam}, hardness amplification \cite{Vio, VW},
probabilistic checkable proofs (PCPs) \cite{ST}, and cryptography
\cite{BOU9, BOU1, CFKLLS, FKS}.

A sum-product problem associated with a graph was initiated by
Erd\H{o}s and Szemer\'{e}di \cite{ES}. Alon et al. \cite{AABL}
studied the sum-product theorems for sparse graphs, and obtained
some nice results when the graph is a matching.

Let us ask does there exist any connection between the
`sum-product problem' and `spectral graph theory'? Surprisingly,
the answer is yes! In fact, the first paper that introduced and
applied the spectral methods to estimate sum-product problems (and
even more general problems) is the paper by Vu \cite{VU} (see also
\cite{HLS}, in which Fourier analytic methods were used to
generalize the results by Vu). In his elegant paper, Vu relates
the sum-product bound to the expansion of certain graphs, and then
via the relation of the spectrum (second eigenvalue) and expansion
one can deduce a rather strong bound. Vinh \cite{VIN} (also see
\cite{VIN2}), using ideas from spectral graph theory, derived a
Szemer\'{e}di-Trotter type theorem in finite fields, and from
there obtained a different proof of Garaev's result \cite{GAR3} on
sum-product estimate for large subsets of finite fields. Also,
Solymosi \cite{SOL2} applied techniques from spectral graph theory
and obtained estimates similar to those of Garaev \cite{GAR} that
already followed via tools from exponential sums and Fourier
analysis. One important ingredient in Solymosi's method
\cite{SOL2} is the well-known {\it Expander Mixing Lemma} (see,
e.g., \cite{HLW}), which roughly speaking, states that on graphs
with good expansion, the edges of the graphs are well-distributed,
and in fact, the number of edges between any two vertex subsets is
about what one would expect for a random graph of that edge
density.

The generalizations of the sum-product problem to polynomials,
elliptic curves, and also the exponentiated versions of the
problem in finite fields were obtained in \cite{ASH, BIW, BT,
CRHA, GASH, SHEN3, SH1, SH2, VU}. Also, the problem in the
commutative integral domain (with characteristic zero) setting was
considered in \cite{VWW}. Some other generalizations to algebraic
division algebras and algebraic number fields were treated in
\cite{BOCH0, CHAN2}. Tao \cite{TAO} settled the sum-product
problem in arbitrary rings.

As we already mentioned, the sum-product theorem is certainly not
true for matrices over $\mathbb{F}_p$. However, Helfgott
\cite{HEL1} proves that the theorem is true for $A \subset
SL_2(\mathbb{F}_p).$ In particular, the set $A\cdot A\cdot A$ is
much larger than $A$ (more precisely, $|A\cdot A\cdot
A|>|A|^{1+\varepsilon}$, where $\varepsilon >0$ is an absolute
constant), unless $A$ is contained in a proper subgroup.
Helfgott's theorem has found several applications, for instance,
in some nonlinear sieving problems \cite{BGS1}, in the spectral
theory of Hecke operators \cite{GJS}, and in constructing
expanders via Cayley graphs \cite{BG1}. Underlying this theorem is
the sum-product theorem. Very recently, Kowalski \cite{KOW}
obtained explicit versions of Helfgott's growth theorem for
$SL_2$. Helfgott \cite{HEL2} proves his result when $A \subset
SL_3(\mathbb{F}_p)$, as well. Gill and Helfgott \cite{GH1}
generalized Helfgott's theorem to $SL_n(\mathbb{F}_p)$, when $A$
is small, that is, $|A|\leq p^{n+1-\delta}$, for some $\delta >0$.
The study of growth inside solvable subgroups of
$GL_r(\mathbb{F}_p)$ is done in \cite{GH2}. Breuillard, Green, and
Tao \cite{BGT}, and also Pyber and Szab\'{o} \cite{PYSZ},
independently and simultaneously, generalized Helfgott's theorem
to $SL_n(\mathbb{F})$, ($n$ arbitrary, $\mathbb{F}$ arbitrary
finite field), and also to some other simple groups, as part of a
more general result for groups of bounded Lie rank; see also
\cite{HRU}.

Let us ask does there exist a `sum-division theory'? Solymosi
\cite{SOL4} using the concept of {\it multiplicative energy}
proved the following: If $A$ is a finite set of positive real
numbers, then
$$|A+A|^{2}|A\cdot A|\geq \frac{|A|^{4}}{4\lceil \log_2 |A| \rceil}.$$
Solymosi's result also gives $$|A+A|^{2}|A/A|\geq
\frac{|A|^{4}}{4\lceil \log_2 |A| \rceil}.$$ Li and Shen \cite{LS}
removed the term $\lceil \log_2 |A| \rceil$ in the denominator. In
fact, they proved the following: If $A$ is a finite set of
positive real numbers, then $$|A+A|^{2}|A/A|\geq
\frac{|A|^{4}}{4},$$ which concludes that
$$\max \{|A+A|,|A/A|\}\geq\frac{|A|^{4/3}}{2}.$$ One may ask about
a `difference-product theory'. The work of Solymosi \cite{SOL1}
considers this type, but the state-of-the-art not only for this
type but also for all combinations of addition, multiplication,
subtraction and division in the case of complex numbers is due to
Rudnev \cite{RUDN1}.

As we already mentioned, Rudnev \cite{RUDN} proved
$$\max \{|A+A|,|A\cdot A|\}\gg \frac{|A|^{12/11}}{\log
^{4/11}|A|},$$ where $A\subset \mathbb{F}_p^*$ with $|A|<\sqrt{p}$
and $p$ large. In Remark 2 of his paper, Rudnev \cite{RUDN},
mentions an interesting fact: ``one can replace either one or both
the product set $A\cdot A$ with the ratio set $A/A$ -- in which
case the logarithmic factor disappears -- and the sumset $A+A$
with the difference set $A-A$"!

The sum-product problem can be applied efficiently to construct
{\it randomness extractors} \cite{BIW, BKSSW, BOU3, BOU7, CG, DGW,
HEHE, YEH}. Inspired by this fact, we are going to discuss some
properties and applications of randomness extractors here. First,
note that all cryptographic protocols and in fact, many problems
that arise in cryptography, algorithm design, distributed
computing, and so on, rely completely on randomness and indeed are
impossible to solve without access to it.

A {\it randomness extractor} is a deterministic polynomial-time
computable algorithm that computes a function \textsf{Ext} :
$\{0,1\}^n \rightarrow \{0,1\}^m$, with the property that for any
defective source of randomness $X$ satisfying minimal assumptions,
\textsf{Ext}($X$) is close to being uniformly distributed. In
other words, a randomness extractor is an algorithm that
transforms a weak random source into an almost uniformly random
source. Randomness extractors are interesting in their own right
as combinatorial objects that ``appear random" in many strong
ways. They fall into the class of ``pseudorandom" objects. {\it
Pseudorandomness} is the theory of efficiently generating objects
that ``appear random" even though they are constructed with little
or no true randomness; see \cite{TREV0, VAD0, VAD1} (and also the
surveys \cite{SHA, SHA1}). Error correcting codes, hardness
amplifiers, epsilon biased sets, pseudorandom generators, expander
graphs, and Ramsey graphs are of other such objects. (Roughly
speaking, an {\it expander} is a highly connected sparse finite
graph, i.e., every subset of its vertices has a large set of
neighbors. Expanders have a great deal of seminal applications in
many disciplines such as computer science and cryptography; see
\cite{HLW, LUBO} for two excellent surveys on this area and its
applications.) Actually, when studying large combinatorial objects
in additive combinatorics, a helpful (and easier) procedure is to
decompose them into a `structured part' and a `pseudorandom part'.

Constructions of randomness extractors have been used to get
constructions of communication networks and good expander graphs
\cite{CRVW, WZ}, error correcting codes \cite{GUR, TZ1},
cryptographic protocols \cite{LU, VAD}, data structures
\cite{MNSW} and samplers \cite{ZUC}. Randomness extractors are
widely used in cryptographic applications (see, e.g., \cite{BTV,
CFPZ, DGKM, DOWI, KLR, KRVZ, LIX1, ZIMA}). This includes
applications in construction of pseudorandom generators from
one-way functions, design of cryptographic functionalities from
noisy and weak sources, construction of key derivation functions,
and extracting many private bits even when the adversary knows all
except $\log^{\Omega(1)}n$ of the $n$ bits \cite{RAO} (see also
\cite{RAO0}). They also have remarkable applications to quantum
cryptography, where photons are used by the randomness extractor
to generate secure random bits \cite{SHA}.

{\it Ramsey graphs} (that is, graphs that have no large clique or
independent set) have strong connections with extractors for two
sources. Using this approach, Barak et al. \cite{BRSW} presented
an explicit Ramsey graph that does not have cliques and
independent sets of size $2^{\log^{o(1)}n}$, and ultimately
beating the Frankl-Wilson construction!

\section*{Acknowledgements}

The author would like to thank Igor Shparlinski for many
invaluable comments and for his inspiration throughout the
preparation of this survey and his unending encouragement. I also
thank Antal Balog, Emmanuel Breuillard, Ernie Croot, Harald
Helfgott, Sergei Konyagin, Liangpan Li, Helger Lipmaa, Devanshu
Pandey, Alain Plagne, L\'{a}szl\'{o} Pyber, Jeffrey Shallit,
Emanuele Viola, and Van Vu for useful comments on this manuscript
and/or sending some papers to me. Finally, I am grateful to the
anonymous referees for their suggestions to improve this paper.

\end{document}